\documentstyle[twoside]{article}
\newtheorem{thm}{Theorem}[section]

\newtheorem{lemma}{Lemma}[section]
\newtheorem{prop}{Proposition}[section]
\newtheorem{cor}{Corollary}[section]
 
\newcommand{\xs}{(x_{1},x_{2},\ldots) }
\begin{document}

\begin{titlepage}
\title{A Simpler Characterization of Sheffer Polynomials}
\author{A. Di Bucchianico\thanks{Author supported by NWO (Netherlands
Organization for Scientific Research).}\\
Department of Mathematics\\
University of Groningen\\
P. O. Box 800\\
9700 AV Groningen, Netherlands\\
A.Di.Bucchianico@math.rug.nl
\and
D. Loeb\thanks{Author supported by URA CNRS 1304.}\\
LABRI\\
Universit\'e de Bordeaux I\\
33405 Talence, France\\
loeb@geocub.greco-prog.fr}
\date{\ }
\end{titlepage}
\maketitle

\begin{abstract}
We characterize the Sheffer sequences by a single convolution identity
$$ F^{(y)} p_{n}(x) = \sum _{k=0}^{n}\ p_{k}(x)\ p_{n-k}(y)$$
where $F^{(y)}$ is a shift-invariant operator. We then study a
generalization of the notion of Sheffer sequences by removing the
requirement that $F^{(y)}$ be shift-invariant. All these solutions can
then be interpreted as cocommutative coalgebras. We also show the connection
with generalized translation operators as introduced by Delsarte. 
Finally, we apply the same convolution to symmetric functions where we
find that the ``Sheffer'' sequences differ from ordinary full divided power
sequences by only a constant factor.
\end{abstract}

\section{Introduction}
The basis of the Umbral Calculus (see \cite{MR} and \cite{RKO}) is the
convolution identity
\begin{equation}\label{conv}
E^{y}p_{n}(x) = \sum _{k=0}^{n}\ p_{k}(x)\ p_{n-k}(y)
\end{equation}
where the {\em shift operator} $E^{y}:K[x]\rightarrow K[x,y]$ is
defined by $E^{y}p(x) = p(x+y)$.  
A {\em sequence of polynomials} is a sequence $(p_{n}(x))_{n=0}^{\infty
}$ of polynomials such that $\deg (p_{n}(x))=n$.  
It is said to be a {\em divided powers sequence} if it
obeys equation \ref{conv}. The Umbral Calculus is the study of such
sequences and their sister sequences of {\em binomial type}
$(q_{n}(x))_{n=0}^{\infty }$ 
with $q_{n}(x) = n!\ p_{n}(x)$ so called since they obey the ``binomial''
identity
$$ E^{y}q_{n}(x) = \sum _{k=0}^{n} {n\choose k}\ q_{k}(x)\ q_{n-k}(y). $$
Famous examples of sequences of binomial type include: the powers of
$x$, the 
lower factorials $x(x-1)\cdots(x-n+1)$, the rising factorials
$x(x+1)\cdots(x+n-1)$, and the Abel polynomials $x(x-na)^{n-1}$.

A related concept is that of the Sheffer sequences. A
Sheffer sequence of polynomials $(s_{n}(x))_{n=0}^{\infty }$ has been
traditionally
defined algebraically by the identity
\begin{equation}\label{shef}
E^{y}s_{n}(x) = \sum_{k=0}^{n}\ p_{n-k}(y)\ s_{k}(x)
\end{equation}
where $(p_{n}(x))_{n=0}^{\infty }$ is itself a divided power sequence
of polynomials.
For example, the Bernoulli polynomials are Sheffer with respect to
$(x^{n}/n!)_{n=0}^{\infty }$.
Thus, {\em a priori}, the Sheffer sequence is a less basic concept than
that of divided power sequences---as far as convolutiion identities is
concerned.

At the Franco-Qu\'ebecois Workshop in May 1991, we asked what sort
of ``shiftless'' Umbral Calculus would arise if the operator $E^{y}$
was replaced by some other shift-invariant operator
$F^{(y)}:K[x]\rightarrow K[x,y]$.
\begin{equation}\label{meq}
F^{(y)}p_{n}(x) = \sum _{k=0}^{n}\ p_{k}(x)\ p_{n-k}(y)
\end{equation}
(We write $F^{(y)}$ so as {\em not} to imply that $F^{(y)}F^{(z)}$ is
necessarily equivalent to $F^{(y+z)}$).
In section \ref{polys}, we will show that only Sheffer sequences obey equation
\ref{meq}. Thus, Sheffer sequences are a much more natural subject of study
than is the special case of divided power sequences. We also show some
connections with the theory of generalized translation operators and Cauchy
problems as presented in \cite{Lev}.

In section \ref{syms}, we seek  parallel results for divided power
sequences of symmetric functions. Surprisingly, up to a constant, the
only ``Sheffer'' sequences of symmetric functions are the divided
power sequences themselves.

We end sections \ref{polys} and \ref{syms} with applications to
coalgebra theory. These results may be safely skipped by any
non-specialist. Solutions to equation~\ref{meq} are interpreted as
cocommutative coalgebras, and classified according to their coalgebraic
properties.

\section{Polynomials}\label{polys}
\subsection{Notation}

Let $K$ be a field of characteristic zero, and let $x,y$ be
indeterminates. Consider a $K$-linear map $\phi: K[x]\rightarrow
K[x]$. Then $\phi $ has a unique $K[y]$-linear extension $\phi '$ to
$K[x,y]$.   By an abuse of notation, we will denote $\phi $ and $\phi
'$ by the same symbol $\phi $. Note that if $\phi $ is a linear map on
$K[x]$ and $\theta  $ is a linear map on $K[y]$, then $\phi $  and
$\theta $ commute when considered as maps on $K[x,y]$. Nevertheless,
$\phi $ and $\theta $ do not necessarily commute with maps $\psi :
K[x]\rightarrow K[x,y]$ such as the shift operator.

Now, consider the natural isomorphism $\pi : K[x]\rightarrow K[y]$.
There is a unique map which we will denote $T_{x}^{y}\phi $ such that
$T_{x}^{y}\phi \circ \pi = \pi \circ \phi $. 
Often a linear map $K[x]\rightarrow K[x]$ will be denoted $\phi _{x}$
with $x$ as a subscript. In that case, $T_{x}^{y}\phi _{x}$ will be
denoted $\phi _{y}$.
For example, if $D_{x}$ is
the derivative with respect to $x$, then $D_{y}$ is the
derivative with respect to $y$. Similarly, if $\epsilon _{x}$ is the 
evaluation map at $x=0$, $\epsilon _{x}p(x)=p(0)$, then $\epsilon
_{y}$ is the evaluation map at $y=0$. Essentially, $\phi_{y} $
behaves with respect to $y$ in the same way as $\phi_{x} $ does with
respect to $x$.

For $n$ a nonnegative integer, let $p_{n}(x)$ be a polynomial of
degree $n$ with coefficients in $K$ and let 
$F^{(y)}$ be a $K$-linear operator $K[x]\rightarrow K[x,y]$ which we
propose as a possible solution to equation \ref{meq}.

Note that we must assume that $F^{(y)}$ is linear in order to 
characterize such operators satisfying equation \ref{meq}, since
equation \ref{meq} specifies the value of
$F^{(y)}$ only  on a basis of $K[x]$.

Finally, we note that in the above notation there are really two kinds
of ``shifts'' $E^{y}:p(x)\mapsto p(x+y)$. If $y$ is taken as a
constant, then $E^{y}:K[x]\rightarrow K[x]$. Whereas, if $y$ is taken
as a variable, then $E^{y}:K[x]\rightarrow K[x,y]$. However,
commutation with one shift guarantees commutation with the other as the next
lemma shows. 

\begin{lemma}
Let $\theta $ be a linear operator on $K[x]$ (and thus on $K[x,y]$).
Then $\theta E^{c}=E^{c}\theta $ for all $c\in K$ if and only if
$\theta E^{y}=E^{y}\theta $.
\end{lemma}

{\em Proof:} {\bf (If)} Trivial, evaluate at $y=c$.

{\bf (Only if)} For any polynomial $p(x),$ the expressions $E^{y}\theta
p(x)$ and $\theta E^{y}p(x)$ agree infinitely often when thought of as
polynomials taking values in $K[x]$. Thus, they must be the same
polynomial, and $\theta E^{y}=E^{y}\theta . \Box $

\subsection{Sheffer Theorem}

In this section, we make the additional assumption that $F^{(y)}$ is
shift-invariant for all $y.$ 

\begin{thm}[Sheffer Theorem] \label{thm1}
Given $F^{(y)}$ and $(p_{n}(x))_{n=0}^{\infty}$ as above, the
following two statements are equivalent. 
\begin{enumerate}
\item $F^{(y)}$ and $(p_{n}(x))_{n=0}^{\infty}$ obey equation \ref{meq}.
\item $P_{x}=\epsilon _{y}\circ F^{(y)}$ is an invertible
shift-invariant operator on $K[x]$,
$(p_{n}(x))_{n=0}^{\infty}$ is  Sheffer relative to the divided power
sequence $(P_{x}^{-1}p_{n}(x))_{n=0}^{\infty}$, and
$F^{(y)}=P_{y} E^{y}.$
\end{enumerate}
\end{thm}

{\em Proof:} {\bf (2 implies 1):}
Define $q_{n}(x)=P_{x}^{-1}p_{n}(x)$ or equivalently $q_{n}(y) =
P_{y}^{-1} p_{n}(y)$. 
By hypothesis, $(q_{n}(x))_{n=0}^{\infty }$ is a divided power
sequence. That is to say, 
$$E^{y}q_{n}(x) = \sum _{k=0}^{n}\ q_{k}(x)\ q_{n-k}(y).$$
Now, operate on both sides of the identity with
$P_{x}P_{y}$. Keeping in mind that $P_{x}E^{y}=E^{y}P_{x}$, we then
get  
$$F^{(y)}p_{n}(x) = P_{y}E^{y}P_{x}q_{n}(x) =  \sum _{k=0}^{n}\
(P_{x}q_{k}(x)) (P_{y}q_{n-k}(y)) = \sum _{k=0}^{n}\ p_{k}(x)\
p_{n-k}(y).$$

{\bf (1 implies 2):} By hypothesis,  $P_{x}$ is shift-invariant, but
we must now show that $P_{x}$ is invertible. Since $p_{0}(0)$ is a
nonzero constant,
$P_{x}p_{n}(x)$ is a polynomial of degree $n$ for all $n$. Hence, $P_{x}$ is
an invertible shift-invariant operator.

We may now let $q_{n}(x)=P_{x}^{-1}p_{n}(x)$ and
$G^{(y)}=P_{y}^{-1} F^{(y)}$. It remains now to show that $G^{(y)}=E^{y}$.

By the above reasoning,
$$ G^{(y)}q_{n}(x)=\sum_{k=0}^{n}\ q_{k}(x)\ q_{n-k}(y). $$
Now, apply $\epsilon _{y}$ to both sides.
Since $\epsilon_{y}G^{(y)} = \epsilon_{y}P_{y}^{-1}P_{y}E^{y}$ is the
identity $I_{x}$, we have
$$ q_{n}(x)=\sum_{k=0}^{n}\ q_{k}(x)\ q_{n-k}(0). $$
In other words, $q_{n}(0)=0$ for $n>0$ and $q_{0}(0)=1$.

Since $F^{(y)}$ and $P_{y}$ are shift-invariant, $G^{(y)}$ is also
shift-invariant. Thus,
we can now apply \cite[Theorem 5.3]{GM} which shows that
$G^{(y)}=E^{y}$. That is to say,  $(q_{n}(x))_{n=0}^{\infty }$ is a divided
power sequence and $(p_{n}(x))_{n=0}^{\infty }$ is Sheffer. $\Box $

\paragraph{}
The above theorem yields interesting Sheffer sequences identities. We
illustrate this with three examples: the Hermite polynomials, the Laguerre
polynomials and the Bernoulli polynomials of the second kind.  Let
$(p_{n}(x))_{n=0}^{\infty}$ be a Sheffer sequence relative to
$(q_{n}(x))_{n=0}^{\infty}$. By the First  Expansion  Theorem
(\cite[Theorem~2]{RKO}),  the  operator $P_{y}$ which maps $q_{n}(y)$ to
$p_{n}(y)$ has expansion $$ P_{y} =
\sum_{k=0}^{\infty}\ p_{n}(0)\ Q_{y}^{n}$$ where $Q_{y}$ is the
delta operator of $(q_{n}(y))_{n=0}^{\infty}$. Note that every shift-invariant
operator can be represented as an integral operator (see \cite{Buc}).

\subsubsection{Hermite}\label{Hermite}
Let $(H_{n}^{\nu}(x))_{n=0}^{\infty}$ be the sequence of Hermite polynomials
of variance $\nu$ where $\nu$ is a real number (see \cite[sect.~10]{RKO}).
Its generating function is
$$ \sum_{k=0}^{\infty}\ H_{k}^{\nu}(x)\ t^{k} = e^{xt-\nu t^{2}/2}.$$
It follows that $P{_y} = e^{-\nu D^{2}/2}$. If $\nu<0$, then
$$ e^{-\nu D^{2}/2}p(x) =\frac{1}{\sqrt{-2 \pi \nu}}\int_{-\infty}^{\infty}
e^{-u^{2}/2\nu}\ p(x+u)\ du,$$
and
$$ \frac{1}{\sqrt{-2 \pi \nu}} \int_{-\infty}^{\infty}
e^{-u^{2}/2\nu} H_{n}^{\nu}(x+y+u)\ du
= \sum_{k=0}^{n}\ H_{k}^{\nu}(x)\ H_{n-k}^{\nu}(y).$$

\subsubsection{Laguerre}
Let  $(L_{n}^{\alpha}(x))_{n=0}^{\infty}$  be  the   sequence   of   Laguerre
polynomials  of  order  $\alpha$  where  $\alpha$  is  a  real  number  (see
\cite[sect.~11]{RKO}).
Its generating function is
$$  \sum_{k=0}^{\infty}\   L_{k}^{\alpha}(x)\   t^{k}   =
(1-t)^{-\alpha-1}\ e^{x\frac{t}{t-1}}.$$
Since  $Q=D/(D-I)$,  it  follows  that  $P{_y}   =   (I-D)^{\alpha+1}$.   If
$\alpha<-1$, then
$$  (I-D)^{\alpha+1}\ p(x)  = \frac{1}{\Gamma(-\alpha-1)}   \int_{0}^{\infty}
t^{-\alpha-2}\ e^{-t}\ p(x+t) dt,$$
and
$$ \frac{1}{\Gamma(-\alpha-1)}  \int_{0}^{\infty}
t^{-\alpha-2}\ e^{-t}\ L_{n}^{\alpha}(x+y+t)\ dt
= \sum_{k=0}^{n}\ L_{k}^{\alpha}(x)\ L_{n-k}^{\alpha}(y).$$

\subsubsection{Bernoulli}
Let  $(b_{n}^{\alpha}(x))_{n=0}^{\infty}$   be   the   sequence   of
Bernoulli polynomials of the second kind. Its generating function is
$$\sum_{k=0}^{\infty}\ b_{k}(x)\ t^{k} =\frac{t}{\log(1+t)}\ (1+t)^{x}.$$
In this case,
$$P_{y}p(x) = \int_{x}^{x+1} p(u) du.$$
Thus, we have
$$ \int_{x}^{x+1} b_{n}(y+u)\ du = \sum_{k=0}^{n}\ b_{k}(x)\ b_{n-k}(y).$$

\subsection{Generalized Sheffer}

Let us now remove the condition that $F^{(y)}$ be shift-invariant
which was so crucial to Theorem \ref{thm1}.
Immediately, we have new solutions to equation \ref{meq}. In fact,
{\em any} sequence of polynomials $(p_{n}(x))_{n=0}^{\infty }$  (with
$\deg p_{n}(x)=n$) gives
rise to a unique operator
of $F^{(y)}$  which verifies equation \ref{meq}.

\begin{thm}[Generalized Sheffer Theorem]\label{gsp}
Let $(p_{n}(x))_{n=0}^{\infty }$ be any sequence of polynomials such that 
$\deg p_{n}(x)=n$. The relation
$$ Qp_{n}(x) = \left\{ \begin{array}{ll}
p_{n-1}(x)              &\mbox{if $n>0$, and}\\
0               &\mbox{if $n=0$}
         \end{array}
\right.
 $$
defines a unique linear operator $Q$. Furthermore, the relations
$$ q_{n}(0) = \delta_{n0},$$
and
$$ Qq_{n}(x) = \left\{ \begin{array}{ll}
q_{n-1}(x)              &\mbox{if $n>0$, and}\\
0               &\mbox{if $n=0$}
         \end{array}
\right.
 $$
define a unique sequence of polynomials $(q_{n}(x))_{n=0}^{\infty }$
which in the philosophy of
\cite{GM} would be called a divided power sequence {\em relative to}
or {\em basic for} $Q$.
The relation
$$ P_{x} q_{n}(x)=p_{n}(x) $$
defines a $K$-linear operator $P_{x}$. (Incidentally, $P_{x}$ is
$Q$-invariant and invertible.)
The only solution $F^{(y)}$ to equation \ref{meq} is $P_{y}G^{(y)}$ where
$G^{(y)}$ is given by the convergent sum
\begin{equation}\label{TayDel}
 G^{(y)} =  \sum _{n=0}^{\infty }\ q_{n}(y)\ Q^{n}.
 \end{equation}
\end{thm}

{\em Proof:} Let us first check that all the objects mentioned above
are well-defined. Since $(p_{n}(x))_{n=0}^{\infty }$ is a basis for
$K[x]$, $Q$ is well defined and  it lowers the
degree of any polynomial by one. Thus, $Q^{-1}$ is well defined {\em
up to a constant}. Since the constant term of $q_{n}(x)$ is given,
$q_{n}(x)$ is well defined. By induction, $Q^{n}$ lowers the degree
of any polynomial by $n$; thus, the sum giving $G^{(y)}$ is in fact
convergent. $P_{x}$ is of course well defined and invertible since
$(q_{n}(x))_{n=0}^{\infty }$ and $(p_{n}(x))_{n=0}^{\infty }$ are both
sequences of polynomials. Thus, $F^{(y)}$ is well defined.

$P_{x}$ is $Q$-invariant because
$$ P_{x}Qq_{n}(x) = P_{x}p_{n}(x) = p_{n-1}(x) = Qq_{n-1}(x) =
QP_{x}q_{n}(x). $$

Again uniqueness of solution is automatic, so it will suffice to
verify that $F^{(y)}$ is in fact a solution.
Now, as in \cite[Lemma 5.2]{GM}, we have
\begin{equation}\label{Gy}
G^{(y)}q_{n}(x) = \sum _{k=0}^{\infty }\ q_{k}(y)\ Q^{k}q_{n}(x) =
\sum _{k=0}^{n}\ q_{k}(y)\ q_{n-k}(x) 
\end{equation}
which given  the $Q$-invariancy of $P_{x}$ can be transformed into
equation \ref{meq} by applying
$P_{x}P_{y}$ to both sides, and exchanging $x$ and $y$. $\Box $

Two explicit examples that illustrate Theorem~\ref{gsp} are:
\begin{itemize}
\item ${\displaystyle p_{n}(x) = 
\frac{(x-1)^{n}}{(n!)^{2}}\ P_{n}~\left( \frac{x+1}{x-1} \right)}$ where
$P_{n}(x)$ is the $nth$ Legendre polynomial.
Here, $Q = D x D = D + x D^{2}$. \cite[Chapter~13, Exercise~11]{Rainville}
\item ${\displaystyle p_{n}(x) = H_{n}(x/2)/(n!)^{2}}$ where $H_{n}(x) =
H_{n}^{1}(x)$ is the $nth$ Hermite polynomial as defined in
Section~\ref{Hermite}. Here, 
$Q =\frac{1}{2} D + \frac{1}{2} x D^{2} - \frac{1}{4} D^{3}$. 
\cite[Chapter~13, p. 220]{Rainville}
\end{itemize}

In these examples, operators of the form $ \sum_{n=0}^{\infty}\
a_{n}(x)\ D^{k}$ 
appear where the $a_{n}(x)$ are polynomials. In fact, any linear operator on
the vector space of polynomials can be represented in this way
(see \cite[Proposition~1]{Kurb}, cf. \cite[Theorems~70 and 77]{Rainville}). The
paper \cite{Kurb} shows an efficient way to calculate the polynomials
$a_{n}(x)$ explicitly.

We see that equation \ref{meq} imposes no conditions
on the sequence $(p_{n}(x))_{n=0}^{\infty }$.
So what does it mean to be a generalized Sheffer sequence if every
sequence is a generalized Sheffer sequence? We can answer this
question as follows. In \cite{GM}, it is shown that any sequence
$(q_{n}(x))_{n=0}^{\infty }$ with $q_{n}(0)=\delta _{n0}$ obeys
equation \ref{meq}.
These sequences are to $F^{(y)}$ as divided power sequences are to the
shift operator. The theorem above says not only that
$(p_{n}(x))_{n=0}^{\infty }$ is
generalized Sheffer, but also how it is so. That is to say, given
$(p_{n}(x))_{n=0}^{\infty }$ there is a unique operator $G^{(y)}$ with
a unique $G^{(y)}$-invariant operator $Q$ with a unique basic sequence
$(q_{n}(x))_{n=0}^{\infty }$. It is this sequence that
$(p_{n}(x))_{n=0}^{\infty }$ is Sheffer with respect to, and $P_{x}$
is the Sheffer operator relating the two sequences (cf. \cite{Sheffer} or
\cite[Chapter~13]{Rainville}). The Hermite example above shows that a suitable
choice of norming constants may change the associated basic sequence. For more
information, see \cite{Bukh90}.

Finally, we note that although equation \ref{meq} does not impose any
conditions on $(p_{n}(x))_{n=0}^{\infty }$, it does impose certain
conditions on the operators $F^{(y)}$.
As seen above, we can multiply $F^{(y)}$ by any invertible operator
$K[x]\rightarrow K[x]$. Moreover, $\epsilon _{y}\circ F^{(y)}$ is
clearly invertible. Thus, without a real loss of generality we can
assume that $\epsilon _{y} \circ F^{(y)}$ is the identity.

Let us partition the set of all linear operators (other than constant multiples
of the identity) according to which operators commute with which operators (cf.
\cite{GM}, \cite{OVV} and \cite[Chapter~13]{Rainville}). 
Then using the methods of \cite{GM} it can be shown that each equivalence class
contains exactly one possible value of $F^{(y)}$ such that $\epsilon_{y}\circ
F^{(y)}$ is the identity.

In particular, the class of shift-invariant operators contains only
such solutions of the form $F^{(y)}=E^{y}$ as we saw above.

\paragraph{}
We now want to point out some connections with generalized translation
operators.  The operators $G^{(y)}$ of equation \ref{TayDel} are generalized
translation operators in the sense of Levitan (see \cite{Lev}). The series on
the right-hand side of equation \ref{TayDel} is called a Taylor-Delsarte series
since they were studied in \cite{Del}.  Levitan stresses the importance of the
infinitesimal generator of the operators $G^{(y)}$.  In our case it is easy to
show that

\begin{prop}
Let $(q_{n}(x))_{x=0}^{\infty }$, $Q$, and $G^{(y)}$ be as in
Theorem~\ref{gsp}. Then we have
$$ \lim_{y \rightarrow 0} \frac{G^{y}-G^{0}}{y} =
\sum_{k=0}^{\infty}\ (D q_{k})(0)\ Q^{k}.$$
\end{prop}

{\em Proof:} Apply the left hand side to the basis
$(q_{n}(x))_{x=0}^{\infty }$. Then equation \ref{meq} yields
$$ \lim_{y \rightarrow 0} \frac{G^{y}-G^{0}}{y}\ q_{n}(x) =
\lim_{y \rightarrow 0}\ \sum _{k=1}^{n}\ q_{n-k}(x)\ q_{k}(y)/y $$
Now, it follows from $q_{k}(0) = \delta_{0k}$ that
$$ \lim_{y \rightarrow 0} \frac{G^{y}-G^{0}}{y}\ q_{n}(x) =
\sum_{k=1}^{n}\ (Dq_{k})(0)\ q_{n-k}(x) 
= \left( \sum _{k=0}^{\infty}\ (Dq_{k})(0)\ Q^{k}\right)q_{n}(x). \Box $$

In particular, it follows from the First Expansion Theorem
\cite[Theorem~2]{RKO} that the right hand side sums to $D$ if $Q$ is a
delta operator with basic set $(q_{n}(x))_{n=0}^{\infty }$.

In \cite{Lev}, Levitan also gives a systematic exposition of the relation
between generalized translation operators and Cauchy problems (i.e, partial
differential equations with initial data). In our case, we have the following
Cauchy problem (cf. \cite[Theorem~5]{Bukh89}):

\begin{prop}[Cauchy problem]\label{Cauchy}
Let $(q_{n}(x))_{x=0}^{\infty }$, $Q$, and $G^{(y)}$ be as in
Theorem~\ref{gsp}, then for all polynomials $p(x)$ we have 
$u(x,y)=G^{(y)}p(x)$ as a solution of the following Cauchy problem
\begin{eqnarray*}
Q_{x}u &=& Q_{y}u\\
u(x,0) &=& p(x)
\end{eqnarray*}
\end{prop}

{\em Proof:} First, note that $u(x,0) = p(x)$ because $G^{(0)} = I$. Since 
$(q_{n}(x))_{n=0}^{\infty }$ is a basis, it suffices to show that 
$Q_{x} G^{(y)} q_{n}(x) =  Q_{y} G^{(y)} q_{n}(x)$. This follows directly from 
equation~\ref{Gy}. $\Box$

\paragraph{}
If $Q=D$ in Proposition~\ref{Cauchy}, then we can easily compute $u$ as
follows: Define new variables $\xi = x$ and $\eta = x+y$. Since $D_{x} =
D_{\xi} + D_{\eta}$ and $D_{y} = D_{\eta}$, the differential equation
transforms into $D_{\xi} u = 0$ with solution $u(x,y) = f(\eta) =
f(x+y)$. Now, 
set $y=0$ which yields $p = f$.  Hence, $u(x,y) = p(x+y)$ as expected.

Another way to solve this Cauchy problem is to proceed as Heaviside did in the
previous century: Fix $x$ and treat $D_x$ as a formal constant. Then the
Cauchy problem becomes an ordinary differential equation whose solution is
readily seen to be $u(x,y) = e^{y\, D_x}\ p(x)$ which equals $p(x+y)$ by the
First Expansion Theorem \cite[Theorem~2]{RKO}.

\paragraph{}
If $Q = c D + x D^{2}$ (cf. the second example below Theorem~\ref{gsp} where $c
= 1$) and $c \geq \frac{1}{2}$, then it follows from
\cite[Theorem~2.4.2.6]{Fein} that 
$$ u(x,y) = \frac{1}{2 \pi} \frac{1}{2 B(c - \frac{1}{2}, \frac{1}{2})}\
\int_{0}^{2 \pi}\ p(x + y - 2 \sqrt{xy}\ \mbox{cos} \phi)\ 
(\mbox{sin}^{2} \phi)^{c-1}\ d\phi$$
where $B$ denotes the beta function.

\noindent
The relation between Cauchy problems and 
generalized translation operators is due to Delsarte
(see \cite{Del}, for recent developments see 
\cite{Markett} and references therein).
Delsarte mainly considered the  Hankel
translation, which is associated with the Sturm-Liouville
operator 
$$ \Delta_{x}= \frac{d^{2}} {dx^{2}} + \frac{2\nu}{x} \frac{d}{dx} .$$
A closed form for the Hankel translation is given by (see e.g.
\cite[p. 4]{Cho})  
$$
G^{(y)}p(x) = \frac{\Gamma(\nu+1/2)} {\Gamma(\nu) \Gamma(1/2)} \int_{0}^{\pi}
p[ \{ x^{2}+y^{2}-2yx\cos\theta\}^{1/2}] (\sin\theta)^{2\nu-1}  d\theta .$$

An Umbral Calculus based on the Hankel translation operator is presented in
\cite{Cho}. This Umbral Calculus is related to Bessel functions.

\subsection{Coalgebra}

The above can be profitably recast in the terminology of coalgebras (see
\cite{Hopf} for the relation between Umbral Calculus and coalgebras).
A {\em coalgebra} is a vector space $V$ equipped with a {\em comultiplication}
$\Delta : V\rightarrow V\otimes V$ and a {\em counitary map} $\epsilon
:V\rightarrow K$. These maps must be coassociative
\begin{equation}\label{coass}
  (I\otimes \Delta )\circ  \Delta =  (\Delta \otimes I)\circ  \Delta
\end{equation}
and obey the counitary property
\begin{equation}\label{counit}
(\epsilon \otimes I)\circ \Delta = I = (I \otimes \epsilon ) \circ \Delta.
\end{equation}

Now, $K[x,y]$ is isomorphic to the tensor product $K[x]\otimes K[x]$,
so any  $F=F^{(y)}$ (satisfying equation~\ref{meq}) would be a
potential candidate for a comultiplication map.
Equation \ref{coass} is automatically satisfied:
\begin{eqnarray*}
(F \otimes I) \circ F p_{n}(x) &=& \sum_{i+j+k=n}\ p_{i}(x)\otimes
p_{j}(x) \otimes p_{k}(x)\\
&=& (I\otimes F) \circ  F p_{n}(x).
\end{eqnarray*}
Moreover,  $F$ is automatically {\em cocommutative} since $\sum
_{k=0}^{n}\ p_{k}(x)\ p_{n-k}(y)$ is symmetric in $x$ and $y$.

By equation \ref{counit},
\begin{eqnarray*}
p_{n}(x) &=& (\epsilon \otimes I)  F p_{n}(x)\\
&=& (\epsilon \otimes I)  \sum _{k=0}^{n}\ p_{k}(x) \otimes p_{n-k}(x)\\
&=& \sum _{k=0}^{n}\ (\epsilon p_{k}(x))\ p_{n-k}(x).
\end{eqnarray*}
Since $\{ p_{n}(x): n\in {\bf N} \}$ is a basis, we have $\epsilon
p_{k}(x) = \delta
_{k0}$. For example, $\epsilon $ is the ``evaluation at zero''
operator if, as in \cite{GM}, $p_{n}(0)=\delta _{n0}$.

We have thus proven the following proposition.
\begin{prop}
All $F^{(y)}$ satisfying equation \ref{meq} define distinct (yet
isomorphic) cohomogeneous cocommutative coalgebras. Conversely, any
cohomogeneous coalgebra isomorphic to $(K[x],E^{y})$ yields a solution
to equation \ref{meq}.
\end{prop}

\begin{cor}
Suppose $F^{(y)}$ together with $(p_{n}(x))_{n=0}^{\infty }$ obeys
equation \ref{meq}, and $F^{(y)}$ together with
$(p'_{n}(x))_{n=0}^{\infty }$ also obeys equation \ref{meq}. Then the
two maps $\epsilon $ and $\epsilon '$ defined by
\begin{eqnarray*}
\epsilon p_{n}(x) &=& \delta _{n0}\\
\epsilon' p'_{n}(x) &=& \delta _{n0}
\end{eqnarray*}
are identical.
\end{cor}

\begin{cor}
Suppose $F^{(y)}$ together with $(p_{n}(x))_{n=0}^{\infty }$ obeys
equation \ref{meq}, and the 
resulting coalgebra is in fact a bialgebra with respect to the usual
multiplication of polynomials. Then $F^{(y)}$ is the map $E^{y-c}$ for
some constant $c$. Thus, $(p_{n}(x))_{n=0}^{\infty }$ is a Sheffer
sequence. The counitary map $\epsilon $ is evaluation at $x=c$. These
bialgebras are then Hopf algebras when equipped with the antipode
$\omega : (y+c)^{n}\mapsto (-1)^{n}(y+c)^{n}. $
\end{cor}

{\em Proof:} If $F^{(y)}$ is an algebra map, then $F^{(y)}$ is the
substitution for $x$ of some polynomial $r(x,y)$. By degree
considerations in equation \ref{meq}, $r(x,y)$ must be of degree one.
Moreover, since $F^{(y)}$ is cocommutative, $r(x,y)$ must be symmetric
in $x$ and $y$. Thus, $r(x,y)=a(x+y)-c$. Consideration of the leading
coefficients in equation \ref{meq} indicates that $a$ must be zero.
Thus, $F^{(y)}=E^{y-c}$. The remaining results are easily
verified. $\Box $

\section{Symmetric Functions}\label{syms}
\subsection{Introduction}
In \cite{L}, the notion (and combinatorial interpretation) of divided
power sequences is extended to the domain of symmetric functions. A {\em
linear divided powers sequence} of symmetric functions
$(p_{n} \xs)_{n=0}^{\infty }$ is a sequence 
of homogeneous symmetric functions---one of each degree---obeying the
following convolution identity
$$ E^{y}p_{n} \xs = \sum _{k=0}^{n}\ p_{k}\xs\ p_{n-k}(y,0,0,\ldots) $$
where the {\em symmetric shift} $E^{y}$ is defined by the rule
$$ E^{y}q \xs = q(y,x_{1},x_{2},\ldots). $$
Well-known examples of linear divided power sequences of symmetric
functions include the elementary $e_{n}\xs $ and complete $h_{n}\xs$
symmetric functions.

Suppose we now generalize to
$$ F^{y}p_{n} \xs = \sum _{k=0}^{n}\ p_{k}\xs\ p_{n-k}(y,0,0,\ldots) $$
where $p_{n}\xs$ is a sequence of homogeneous symmetric functions---one
for each degree---and  $F^{y}$ is a linear operator. In this case,
there is not much to say about $F^{y}$. It is not defined on a
basis, so there are not enough constraints to characterize it
completely.

Clearly, we are considering the wrong generalization of polynomial
sequences. We must turn to the subject of \cite{L2}, {\em full
sequences of symmetric functions}, since it is those sequences which
serve as a useful basis for the space of symmetric functions.

\subsection{Notation}

A {\em partition} $\lambda $ is an eventually zero, decreasing
sequence of natural numbers $\lambda _{1}\geq \lambda _{2}\geq \cdots
=0$.  Its conjugate, denoted $\lambda '$, is defined by the rule
$$ \lambda' _{i} = \left| \{j: \lambda _{j}\geq i \}\right| .$$

We will compare partitions and/or vectors in two different ways.
\begin{itemize}
\item First, they can be compared coordinate wise: $\alpha \leq \beta $ if
and only if $\alpha_{i}\leq \beta _{i}$ for all $i$.
\item Second, they can be compared using the reverse lexicographical
order. That is to say, they are ordered as if they were words written
in Hebrew or Arabic (from right to left). $\alpha\ll \beta $ if and
only if there is an $i$ such that $\alpha_{i}<\beta _{i}$ and
$\alpha_{j}=\beta _{j}$ for all $j>i$.
\end{itemize}

Let ${\cal P}$ be the set of all partitions and ${\cal P}_{n}$ be the
set of all partitions summing to $n$.
Clearly, only $\ll $ is a total ordering of 
${\cal P}$. In fact, $\ll $ is a strengthening of the $<$ relation
which itself is so weak as to be equality when restricted to  ${\cal
P}_{n}$. 

The {\em monomial symmetric functions} $m_{\lambda }\xs $ for $\lambda
\in {\cal P}_{n}$ form a basis for the vector space of homogeneous
symmetric functions of degree $n$. In fact, $(m_{\lambda '}\xs
)_{\lambda \in{\cal P}}$ will be
our canonical example of a full sequence (just as
$(x^{n})_{n=0}^{\infty }$ is the
typical sequence of polynomials).

In general, in a {\em full sequence}
$(p_{\lambda }\xs )_{\lambda \in{\cal P}}$, the symmetric functions
$p_{\lambda }\xs $  must be homogeneous of degree $n$ (for
$\lambda \in {\cal P}_{n}$). Moreover, they must have
expansions in terms of the monomial symmetric functions
whose index follows $\lambda '$ in reverse lexicographical order
\begin{equation}\label{full}
p_{\lambda}\xs =\sum _{\mu \underline{\gg}\ \lambda '} b_{\lambda \mu
} m_{\mu }\xs
\end{equation}
where $b_{\lambda \lambda }$ is never zero.

A full sequence is thus a basis for the space of symmetric functions.

Even though $p_{\lambda} \xs $ is only defined for $\lambda $ a
partition, it will be convenient to extend its definition to all
vectors of integers with finite support. If $\alpha_{i}$ is always
nonnegative, then there is a unique partition $\lambda $ which is a
permutation of $\alpha$. We then write
$$ p_{\alpha}\xs =p_{\lambda }\xs . $$
On the other hand, if $\alpha_{i}<0$ for some $i$, we write
$$ p_{\alpha}\xs =0. $$

Finally, we must define a few linear operators; the multivariate symmetric
derivative $D_{\lambda }$ is most simply defined by
$$ D_{\lambda }m_{\mu }\xs = m_{\mu -\lambda } $$
while the augmentation $\epsilon $ is defined by
$$ \epsilon p\xs = p(0,0,\ldots). $$
Note that $E^{a}=\sum_{n=0}^{\infty } a^{n}D_{(n)}.$
A linear operator $\theta $ is said to be shift-invariant is
$E^{a}\theta =\theta E^{a}$. In that case, we have the following
convergent expansion of $\theta $ in terms of $D_{\lambda }:$ 
$$ \theta =\sum _{\lambda }\ \epsilon (\theta m_{\lambda }\xs )\ 
D_{\lambda }. $$

Now, we can define the object of interest; a {\em full divided powers
sequence} is a full sequence of symmetric functions $(p_{\lambda
}\xs)_{\lambda \in {\cal P}}$ which obeys the convolution identity
$$ E^{y}p_{\lambda } \xs = \sum _{\alpha}\ p_{\alpha}\xs\ p_{\lambda
-\alpha}(y,0,0,\ldots) $$ 
where the sum is over all integer vectors $\alpha$ with finite
support. 

\subsection{Sheffer Theorem}

What linear operators $F^{y}$ and full sequences $p_{\lambda }\xs $
obey
\begin{equation}\label{seq}
F^{y}p_{\lambda } \xs = \sum _{\alpha}\ p_{\alpha}\xs\ p_{\lambda
-\alpha}(y,0,0,\ldots)?
\end{equation}

\begin{thm}[Sheffer Theorem]
Given that $F^{y}$ is a shift-invariant operator obeying equation
\ref{seq}, then $F^{y}=cE^{y}$. In other words, $(p_{\lambda
}\xs)_{\lambda \in {\cal P}} $ is
up to constant factor equal to a full divided power
sequence of symmetric functions. 
\end{thm}

{\em Proof:} First, consider $F^{0}$.
$$F^{0}p_{\lambda } \xs = \sum _{\alpha}\ p_{\alpha}\xs\ p_{\lambda
-\alpha}(0,0,0,\ldots).$$
However, for $\alpha\neq (0)$, $p_{\alpha}(0,0,\ldots)$ is zero while
$p_{(0)}(0,0,\ldots)=c\neq 0.$

Without loss of generality, we can assume that $c=1$. Otherwise,
replace $F^{y}$ with $\frac{1}{c}F^{y}$ and $p_{\lambda }\xs $ with
$\frac{1}{c}p_{\lambda }\xs $. It remains then to show that $F^{y}=E^{y}$.

Since $F^{y}$ and $E^{y}$ are both shift-invariant, so is their
difference which we can then expand in the form
$$ F^{y}-E^{y} = \sum _{\lambda }\ c_{\lambda }D_{\lambda }. $$
We will show by induction on $\lambda $ (ordered reverse
lexicographically) that $c_{\lambda }=0$ and thus $F^{y}=E^{y}$. 
The base case $\lambda =(0)$ has already been dispensed with.

Let $\lambda \in {\cal P}_{n}$ $(n>0)$, and suppose that $c_{\mu }=0$
for $\mu \ll \lambda \in {\cal P}_{n}$ and for $\mu \in {\cal P}_{m}$ with
$m<n$.  We must show that $c_{\lambda }=0$.
By induction,
$$ (F^{y}-E^{y})p_{\lambda '}\xs =c_{\lambda }b_{\lambda '\lambda '} $$
where the $b$ sequence is defined by equation \ref{full}.
However, the right hand side is equal to
$$ p_{\lambda '}\xs -p_{\lambda '}(y,x_{1},x_{2},\ldots) + \sum
_{\alpha\neq (0)}\ p_{\lambda '-\alpha}\xs\ p_{\alpha}(y,0,0,\ldots) $$
which is homogeneous of degree $n$ in the variables
$x_{1},x_{2},\ldots,$ and $y$. Thus, the right hand side has no
constant term. Therefore, the constant $c_{\lambda }b_{\lambda
'\lambda '}$ must be zero. However, $b_{\nu \nu }$ is never zero, so
we must have $c_{\lambda }=0. \Box $

{\em Open Problem:} What happens if we no longer assume that $F^{y}$
is shift-invariant? Do we get an analog of Proposition \ref{gsp} ?

\subsection{Coalgebra}
As seen in \cite{L}, all operators of the form $F^{y}$ obeying
\ref{seq} serve as the comultiplication of a (stronly) cohomogeneous
cocommutative Hopf 
algebra over the symmetric functions, and conversely. For the
symmetric shift operator, for example, the augmentation $\epsilon $ is
the counitary map, and the antipode is the classical involution of
symmetric functions
$$ \omega h_{n}\xs = (-1)^{n}e_{n}\xs .  $$


\begin{thebibliography}{99}
\bibitem{Buc} {\sc A. Di Bucchianico}, {\em Representations of Sheffer
polynomials}, submitted to J. Math. Anal. Appl.
\bibitem{Bukh89} {\sc V. M. Bukhstaber and A. N. Kholodov}, {\em
Groups of formal 
diffeomorphisms of the superline, generating functions for sequences of
polynomials, and functional equations}, Izv. Akad. Nauk SSSR {\bf 53} (1989),
944-970 (English transl. in Math. USSR Izvestiya {\bf 35} (1990), 277-305; 
MR~91h:58014).
\bibitem{Bukh90} {\sc V. M. Bukhstaber and A. N. Kholodov}, {\em Boas-Buck
structures on sequences of polynomials}, Funct. Anal. Appl. {\bf 23 (4)}
(1990), 266-276 (MR~91d:26017).
\bibitem{Cho} {\sc F. M. Cholewinski}, {The Finite Calculus associated
with Bessel Functions}, Contemporary Mathematics {\bf 75}, American
Mathematical Society, 1988.
\bibitem{Del} {\sc J. Delsarte}, {\em Sur une extension de la formule de
Taylor}, J.  Math. Pur. Appl. {\bf 17} (1938), 213-231.
\bibitem{Fein} {\sc P. Feinsilver and R.~Schott}, {\em Algebraic structures
and operator calculus}, to appear: Kluwer, The Netherlands.
\bibitem{Kurb} {\sc S. G. Kurbanov and V. M. Maksimov}, {\em Mutual
expansions of 
differential operators and divided difference operators}, Dokl. Akad. Nauk
UzSSR {\bf 4} (1986), 8-9 (MR~87k:05021).
\bibitem{Lev} {\sc B.M. Levitan}, {\em Generalized Translation Operators},
Israel Program Sci. Translations, Jerusalem, 1964.
\bibitem{L} {\sc D. Loeb}, {\em Sequences of symmetric functions of binomial
type}, Stud. Appl. Math. {\bf 83} (1990), 1-30.
\bibitem{L2} {\sc D. Loeb}, {\em Sequences of symmetric functions
of binomial type II: Full Sequences}, In progress.
\bibitem{Markett} {\sc C. Markett}, {\em A new proof
of Watson's product formula for Laguerre polynomials
via a Cauchy problem associated with a singular
differential operator}, SIAM J. Math. Anal. {\bf 17}
(1986), 1010-1032.
\bibitem{GM} {\sc G. Markowsky}, {\em Differential operators and the theory of
binomial enumeration}, J. Math. Anal. Appl. {\bf 63} (1978), 145-155.
\bibitem{IGM} {\sc I. G. Macdonald}, {\em Symmetric Functions and Hall
Polynomials}, Oxford Mathematical Monographs, Clarendon Press, Oxford,
1979.
\bibitem{Hopf} {\sc R. Morris (ed.)}, {\em Umbral Calculus and Hopf
Algebras}, Contemporary Mathematics {\bf 6}, American Mathematical Society,
1982.
\bibitem{MR} {\sc R. Mullin and G.-C. Rota}, {\em On the Foundations of
Combinatorial Theory: III. Theory of Binomial Enumeration}, in: {\sc B.~Harris
(ed.)}, {\em Graph Theory and Its Applications}, Academic Press, 1970, 167-213.
\bibitem{Rainville} {\sc E.D. Rainville}, {\em Special Functions}, MacMillan,
New York, 1960.
\bibitem{RKO} {\sc G.-C. Rota, D. Kahaner and A. Odlyzko}, {\em On the
Foundations of Combinatorial Theory: VIII. Finite Operator Calculus}, J. Math.
Anal. Appl. {\bf42} (1973),  684-760.
\bibitem{Sheffer} {\sc I. M. Sheffer}, {\em Some properties of polynomial sets
of type zero}, Duke Math. J. {\bf 5} (1939), 590-622 (MR~1, 15). 
\bibitem{OVV} {\sc O. V. Viskov}, {\em Operator characterization of
generalized Appell polynomials}, Dokl. Akad. Nauk SSR {\bf 225} (1975), 749-752
(English transl. in Soviet Math. Dokl. {\bf 16} (1975), 1521-1524).
\end{thebibliography}
\end{document}